\documentclass[12pt,a4paper]{amsart}

\usepackage[left=3cm,right=3cm,top=3.5cm,bottom=3.5cm]{geometry}

\usepackage{amsmath,amsthm,amssymb}
\usepackage[english]{babel}

\newtheorem{theorem}{Theorem}[section]

\newtheorem{prop}[theorem]{Proposition}

\newtheorem{coro}[theorem]{Corollary}

\theoremstyle{definition}
\newtheorem{dfn}[theorem]{Definition}
\newtheorem{rem}[theorem]{Remark}
\newtheorem{exam}[theorem]{Example}

\numberwithin{equation}{section}

\DeclareMathOperator{\ord}{ord}

\newcommand{\rleft}{\mathopen{}\mathclose\bgroup\left}
\newcommand{\rright}{\aftergroup\egroup\right}

\def\T{{\Bbb T}}
\def\C{{\Bbb C}}

\def\Q{{\Bbb Q}}
\def\R{{\Bbb R}}
\def\P{{\Bbb P}}
\def\Z{{\Bbb Z}}

\begin{document}

\title[On the stringy Hodge numbers of 
 mirrors]
{On the stringy Hodge numbers of 
 mirrors of quasi-smooth Calabi-Yau hypersurfaces}

\author{Victor V. Batyrev}
\address{Fachbereich Mathematik, Universit\"at T\"ubingen, Auf der
Morgenstelle 10, 72076 T\"ubingen, Germany}
\curraddr{}
\email{batyrev@math.uni-tuebingen.de}
\thanks{}

\thanks{}
%\date{\today}
\begin{abstract}
Mirrors $X^{\vee}$ of quasi-smooth Calabi-Yau hypersurfaces 
$X$ in weighted projective spaces 
$\P(w_0, \ldots, w_d)$ can be obtained as Calabi-Yau compactifications of non-degenerate affine toric hypersurfaces  defined by 
Laurent polynomials whose Newton polytope is 
the lattice simplex 
spanned by $d+1$ lattice vectors $v_i$ satisfying the 
relation $\sum_i w_i v_i =0$. In this paper, we 
compute the stringy $E$-function of mirrors $X^\vee$ and 
compare it with the Vafa's orbifold $E$-function of quasi-smooth Calabi-Yau hypersurfaces $X$. 
As a result, 
we prove 
the equalities of Hodge numbers $h^{p,q}_{\rm str}(X^{\vee}) = h^{d-1-p,q}_{\rm orb}(X)$ 
 for all $p, q$ and $d$ as it is
expected  in mirror symmetry. 

 \end{abstract}

\maketitle

\section{Introduction}

Let $X$ be a $d$-dimensional Calabi-Yau variety, i.e., a normal 
irreducible projective variety with at worst 
canonical Gorenstein singularities satisfying the conditions: $\Omega_X^d \cong {\mathcal O}_X$ and $h^i(X, {\mathcal O}_X) = 0 $ $(0 < i < d)$. We take a desingularization   $\rho\, :\, Y \to X$ 
of $X$ such that the exceptional 
locus of $\rho$ is a union of smooth irreducible divisors
$D_1, \ldots, D_s$ with only normal crossings and write 
\[ K_Y = \rho^* K_X + \sum_{i=1}^s a_i D_i \]
with non-negative integers $a_i$. We set $D_{\emptyset}:= Y$ and consider for any subset 
$J \subseteq I = \{1, \ldots, s \}$ the  
smooth projective  
variety $D_J: = \bigcap_{j \in J} D_j \subseteq Y$ together 
with its open subvariety 
$D_J^{\circ} := D_J \setminus \bigcup_{j \not\in J} D_j$. We denote 
by $E(D_J; u,v):= \sum_{p,q} (-1)^{p+q}h^{p,q}(D_J) u^pv^q$ 
the usual Poincar{\'e} polynomial of Hodge numbers of 
smooth projective varieties $D_J$ and extend them by additivity 
to $E$-polynomials $E(D_J^\circ; u, v):=  \sum_{p,q} 
e^{p,q}(D^\circ_J) u^pv^q$ defined by the Hodge-Deligne numbers 
$e^{p,q}(D^\circ_J)$ of quasi-projective varieties $D_J^\circ$. 

\begin{dfn} \cite{Bat98}
The stringy $E$-function of a Calabi-Yau variety $X$ is  
the rational function in two variabels $u, v$ defined 
by two equivalent formulas: 
\begin{align*}
E_{\rm str}(X; u,v)  & :=  \sum_{\emptyset \subseteq J \subseteq I}
\left( \prod_{j \in J} \frac{uv-1}{(uv)^{a_j+1} - 1} -1 \right)
E(D_J; u,v)  \\ 
& =  \sum_{\emptyset \subseteq J \subseteq I}
\left( \prod_{j \in J} \frac{uv-1}{(uv)^{a_j+1} - 1} \right)
E(D_J^\circ; u,v).
\end{align*}
  If $E_{\rm str}(X; u,v)$ is 
a polynomial, then the {\em stringy Hodge numbers} $h^{p,q}_{\rm str}(X)$ are defined via its coefficients:
\[ E_{\rm str}(X; u,v)= \sum_{p,q}    (-1)^{p+q}h_{\rm str}^{p,q}(X) u^pv^q. \] 

\end{dfn}

The stringy Hodge numbers are very 
useful for testing Mirror Symmetry. If two $d$-dimensional Calabi-Yau varieties $X$ and $X^\vee$ 
are mirror symmetric to each other, then 
\begin{align} \label{pq}
h^{p,q}_{\rm str}(X) =  h^{d- p,q}_{\rm str}(X^\vee), 
\;\; \forall p,q.
\end{align} 
Equations (\ref{pq}) can be equivalently reformulated 
 in the form 
 \begin{align} \label{E-pq}
E_{\rm str}(X^\vee; u,v) = (-u)^d E_{\rm str}(X; u^{-1}, v)  .
\end{align} 
Equation (\ref{E-pq}) has been checked for all 
pairs of Calabi-Yau hypersurfaces
in Gorenstein toric Fano varieties defined by a pair $(\Delta, \Delta^*)$ 
of dual to each other reflexive polytopes \cite{Bat94,BB96}. 
Note that in general the stringy $E$-function 
 $E_{\rm str}(X; u,v)$ is a rational 
function whose denominator is a product of cyclotomic polynomials $\Phi_k(uv)$ in $uv$. 
Therefore, the equation (\ref{E-pq}) 
can be satisfied only provided   $E_{\rm str}(X; u,v)$ is a polynomial.

In this paper we consider Calabi-Yau hypersurfaces  $X= X_w$  
 of degree $w = \sum_{i=0}^d w_i$ 
in weighted projective space $\P(\overline{w})= 
\P(w_0, w_1, \ldots, w_d)$ which is not assumed 
to be  
Gorenstein.  Such a weighted projective 
 space $\P(\overline{w})$ is defined by  a {\em weight vector} 
 $\overline{w} \in \Z^{d+1}$ having  {\em $IP$-property}.
 Recall that  $\overline{w} \in \Z^{d+1}$ has $IP$-property
 if     
the convex hull 
of all non-negative lattice points $(u_0, u_1, \ldots, u_d) \in 
\Z_{\geq 0}^{d+1}$ satisfying the relation 
$\sum_{i=0}^d w_i u_i =w$   is a $d$-dimensional lattice 
polytope  
containing ${\bf 1} := (1, \ldots, 1)$ in its interior \cite{Ska96}.  
This lattice polytope  is just the 
Newton polytope 
of a general weighted homogeneous polynomial 
$W$ defining $X_w$ in $\P(\overline{w})$. 
If $0 \in \C^{d+1}$ is the only singularity of the affine hypersurface $\{W = 0 \} \subset \C^{d+1}$, then the corresponding 
Calabi-Yau hypersurface $X_w \subset \P(\overline{w})$ is called {\em 
quasi-smooth} and the weight vector $\overline{w}$ is called 
{\em transverse}. It is known that any 
transverse weight vector $\overline{w}$ has $IP$-property 
\cite{Ska96}. 
 We consider the weighted 
projective space $\P(\overline{w})$ as a toric variety
with respect to the natural regular action of the $d$-dimensional algebraic 
torus $\T_{\overline{w}}:=(\C^*)^{d+1}/H$, where 
$$H =\{ (\lambda^{w_0}, \ldots, \lambda^{w_d}) \in (\C^*)^{d+1}\; : \;  \lambda \in \C^*\}.$$ 
The lattice of characters of $\T_{\overline{w}}$ is the sublattice 
\[ N_{\overline{w}}:= \{ (u_0,u_1, \ldots, u_d) \in \Z^{d+1} \, |\, \sum_{i=0}^d w_i u_i = 0 \} \subset \Z^{d+1}. \] 
We use the following common notations:
\[ (q_0, q_1, \ldots, q_d) := 
\left(\frac{w_0}{w}, \frac{w_1}{w}, \ldots, \frac{w_d}{w} \right) \in \Q_{>0}^{d+1}, \]
\[ \theta(l):=(\theta_0(l), \theta_1(l), \ldots, \theta_d(l)) = 
 (lq_0, lq_1, \ldots, lq_d), \; \; 0 \leq l < w; \] 
\[ \widetilde{\theta}(l):= (\widetilde{\theta}_0(l), 
\widetilde{\theta}_1(l), \ldots, \widetilde{\theta}_d(l)) \in [0,1)^{d+1}, \; \; 0 \leq l < w, \]
where $\widetilde{\theta}(l)$ is the canonical representative 
of $\theta(l)$ in $(\R/\Z)^{d+1}$. 

Let  $\T_{\overline{w}}^* \cong (\C^*)^d$ be the 
 dual to $\T_{\overline{w}}$ torus with the  lattice of characters $M_{\overline{w}} = 
 \Z^{d+1}/\Z\overline{w}$. Denote by  
 $\Delta_{\overline{w}}^*$ the $d$-dimensional lattice simplex 
 in $M_{\overline{w}} \otimes \R$ whose vertices  
  $v_0, v_1, \ldots, v_d$ are primitive lattice vectors spanning   
$M_{\overline{w}}$ and satisfying the relation $\sum_{i=0}^d w_i v_i = 0$. 
We claim that mirrors  
$X^\vee_{\overline{w}}$ of quasi-smooth Calabi-Yau hypersurfaces $X_w \subset
\P(\overline{w})$  can be obtained as Calabi-Yau  
compactifications of non-degenerate affine hypersurfaces 
in $Z_{\overline{w}} \subset \T_{\overline{w}}^*$ 
defined by Laurent polynomials with the  Newton polytope 
$\Delta_{\overline{w}}^*$.
Using 
a formula of Vafa for the orbifold Euler 
number of a quasi-smooth Calabi-Yau hypersurface $X_w \subset 
\P(\overline{w})$ \cite{Vaf89}: 
\[  \chi_{\rm orb}(X_w) =  
 \frac{1}{w} \sum_{l,r=0}^{w-1} \, \,
\prod_{0 \leq i \leq d \atop  lq_i,rq_i \in \Z}
\left( 1 - \frac{1}{q_i} \right), \]
the topological mirror symmetry 
duality  of  the stringy 
Euler numbers  
\[  E_{\rm str }(X^\vee_{\overline{w}}; 1,1)  =(-1)^{d-1}  
E_{\rm str}(X_w; 1,1)  \]
was checked for this mirror 
construction  in  \cite[Theorem 1.6]{BS20}.  
We note that the quasi-smooth Calabi-Yau hypersurfaces $X_w \subset 
\P(\overline{w})$ are orbifolds and, by a theorem of Yasuda \cite{Yas04}, the stringy 
Hodge numbers $h^{p,q}_{\rm str}(X_w)$ are equal to 
the orbifold Hodge numbers $h^{p,q}_{\rm orb}(X_w)$ 
introduced by Chen and Ruan \cite{CR04}. 

In order 
to compute the orbifold Hodge numbers $h^{p,q}_{\rm orb}(X_w)$ 
by the method of Vafa 
and to compare them with the stringy Hodge numbers 
of mirrors $X_{\overline{w}}^\vee$ we need 
some integers associated with 
elements $l$ in the cyclic  group $G_0:= \Z/w\Z$: 
\begin{itemize}
\item ${\rm age}(l) := \sum_{i =0}^{d} \widetilde{\theta}_i(l) = 
\sum_{\widetilde{\theta}_i(l) \neq 0}  \widetilde{\theta}_i(l), \;\; (l \in \Z/w\Z)$;   
\item ${\rm size}(l) := 
{\rm age}(l) + {\rm age}(w - l) = {\rm Card}\{ i \in \{0,1, \ldots, d\} \, :\,  
\widetilde{\theta}_i(l) \neq 0  \}$. 
\end{itemize} 

If $P(t) = \sum_i b_i t^{r_i}$ is a polynomial or a power series with exponents $r_i \in \frac{1}{w}\Z$ we denote by $[ P(t)]_{int}$ the projector dropping all monomials 
with fractional exponents in $P(t)$
and leaving  only monomials $b_it^{r_i}$ with exponents 
$r_i \in \Z$.  Using 
the action of $G_0$ by $w$-th root of unity 
$e^{2\pi ik/w}$ on $t^{k/w}$, the projector  $[ *]_{int}$ can be expressed 
as a standard Reynolds operator 
\[ [ P]_{int} = \frac{1}{|G_0|} \sum_{g \in G_0} P^g.  \]
In particular, one has  
\begin{align} \label{proj}
[ P \cdot Q ]_{int} = P \cdot [Q]_{int},  
\end{align} 
if $P$ contains only monomials with integral exponents. 

Our main result that extends  \cite[Theorem 1.6]{BS20} is 
the following:

\begin{theorem} \label{main}
Let $\overline{w} \in \Z^{d+1}$ be a weight vector 
with $IP$-property. Denote by $Z_{\overline{w}} \subset \T^*_{\overline{w}}$ 
a  non-degenerate 
affine hypersurface defined 
by a Laurent polynomial with the Newton 
polytope $\Delta_{\overline{w}}^* = {\rm conv}(v_0, \ldots, v_d)
$. Then $Z_{\overline{w}}$ admits 
a Calabi-Yau compactification $X^\vee_{\overline{w}}$ and the 
stringy $E$-function of $X^\vee_{\overline{w}}$ can be computed by the formula
\[ E_{\rm str}(X^\vee_{\overline{w}}; u,v) = 
\frac{1}{uv}  \sum_{0 \leq l < w}  
 \left[ \prod_{\widetilde{\theta}_i(l) = 0} \frac{(uv)^{q_i}  - uv}{1 - (uv)^{q_i}}  \right]_{int}
\cdot   (- u)^{{\rm size}(l)} 
\left( \frac{v}{u} \right)^{{\rm age}(l)}  .\]
Moreover, if the weight vector 
$\overline{w}$ is transverse, then  $E(X^\vee_{\overline{w}}; u,v)$ is 
a polynomial and one has
\[ E_{\rm str}(X^\vee_{\overline{w}}; u,v) = (-u)^{d-1}E_{\rm orb}(X_w; u^{-1},v) \]
where $X_w \subset \P(\overline{w})$ is a quasi-smooth 
Calabi-Yau hypersurface of degree $w = \sum_{i=0}^d w_i$.  
\end{theorem}

We remark that the considered  mirror construction for arbitrary 
quasi-smooth Calabi-Yau hypersurfaces $X_w \subset \P(\overline{w})$ extends the well-known  Berglund-H\"ubsch-Krawitz mirror construction \cite{BH93,Kra10,CR18}. In this construction,
one considers a transverse polynomial $W$ of special form 
\[ W = \sum_{i=0}^d c_i \prod_{j =0}^d z_i^{a_{ij}}, \;\; c_0, c_1, \ldots, c_d \in \C^*. \]
Such a polynomial $W$ is called {\em invertible}.  Consider 
the action of  $l \in G_0 = \Z/w\Z$  
on $\C^{d+1}$ by diagonal matrices 
\[ {\rm diag}(e^{2\pi i lq_0}, e^{2\pi i lq_1},\ldots, e^{2\pi i lq_d}). \]
One associates with the quasi-smooth Calabi-Yau hypersurface $X_w = \{ W = 0 \} \subset \P^d(\overline{w})$  an orbifold Landau-Ginzburg  model  
with the $G_0$-invariant superpotential $W$.   
We set  $\widetilde{G_0} := SL(d+1, \C) \cap 
\widetilde{G}_{\widetilde{W}}$, where  $\widetilde{G}_{\widetilde{W}} \cong {G}_{{W}}  $ is a maximal abelian 
diagonal symmetry group of the dual invertible polynomial
\[ \widetilde{W} := \sum_{j =0}^d c_j \prod_{i=0}^d 
\widetilde{z}_i^{a_{ij}}. \]
Since every monomial $W_i$ in $W$ has degree $w$, we the obtain the equations  
\[  \sum_{j=0}^d a_{ij}w_j = w \;\; \forall  i \in \{0,1, \ldots, d\} \]
which imply that the Newton polytope of the Laurent polynomial 
\[ \widetilde{W}/ \left( \prod_{i=0}^d \widetilde{z}_i \right)^w \]
with respect to the 
$\widetilde{G_0}$-invariant lattice $M_{\widetilde{w}} = \Z^{d+1}/\Z\widetilde{w}$  
is a $d$-dimensional simplex $\Delta_{\overline{w}}^*$ 
whose vertices $v_0, v_1, \ldots, v_d$ generate 
$ M_{\widehat{w}}$ and satisfy the relation $\sum_{i=0}^d w_i v_i =0$. 

\begin{rem}
Note that the mirror construction for quasi-smooth hypersurfaces 
$X_w \subset \P(\overline{w})$ and the proof of Theorem \ref{main} naturally extend to the case of  Calabi-Yau quotients
$X_w/G \subset \P(\overline{w})/G$, where $G \subset 
SL(d+1, \C)$ $(G_0 \subset G)$ is a
finite abelian group of diagonal symmetries 
of the defining transverse polynomial $W$. In this case, 
mirrors of $X_w/G$ are Calabi-Yau compactifications of 
non-degenerate affine hypersurfaces $Z_w \subset \T^*_{\overline{w}}(G)$ whose Newton polytope 
is the same lattice simplex $\Delta_{\overline{w}}^*$, but considered 
with respect to a larger lattice $M_{\overline{w}}(G)$ such that $M_{\overline{w}}(G)/M_{\overline{w}}$ is 
isomorphic to $G/G_0$. 
\end{rem}

\begin{rem}
A combinatorial formula for the orbifold $E$-function  
of pair $(W, G)$, in spirit of Vafa \cite{Vaf89}, has been obtained by 
Ebeling and Takahashi in \cite[Theorem 6]{ET13}. The topological mirror duality between the  orbifold $E$-functions of $(W, G)$ and $(\widetilde{W}, \widetilde{G})$ was shown by Ebeling, Gusein-Zade, and Takahashi in  \cite[Theorem 9]{EGZT16}. 
The isomorphism of the orbifold Hodge spaces in Berglund-H\"ubsch-Krawitz mirror construction  was shown by Chiodo-Ruan \cite{CR10} and Borisov \cite{Bo13}.  
\end{rem}

\section{$E$-polynomials of non-degenerate hypersufaces}

In this section we review some results of Danilov and 
Khovanski\v{i} \cite{DKh86} (see also \cite{Bat93}). 

Let $X$ be a  quasi-projective variety over $\C$. 
Then the cohomology space with compact supports 
$H^k_c(X, \C)$ have natural 
mixed Hodge structure  and one defines {\em Hodge-Deligne numbers} of $X$ as 
\[ e^{p,q}(X) := \sum_{k} (-1)^k h^{p,q}(H_c^k(X)). \]
Then 
\[ E(X; u,v):= \sum_{p,q} e^{p,q}(X) u^pv^p . \]
is called {\em $E$-polynomial} of $X$. 

Let $M$ be a free abelian group of rank $d$. We consider 
$M$ as the lattice of characters of $d$-dimensional 
algebraic torus $\T:= {\rm Hom}(M, \C^*)$, i.e., $\T = 
{\rm Spec}\,  \C[M]$. A Laurent polynomial $f \in \C[M]$ 
defines an affine hypersurface $Z_f = \{ f = 0\} \subset \T$. 
Assume that 
the Newton polytope  $\Delta$ of $f$ has dimension $d$. 
By choosing a $\Z$-basis of $M$, we obtain an isomorphism 
$\C[M] \cong \C[t_1^{\pm 1}, \ldots, t_d^{\pm 1}]$. 
Denote by $M_\Delta:= C_\Delta \cap (\Z \oplus M)$ the graded monoid of all 
lattice  points in   the $(d+1)$-dimensional cone 
$C_\Delta := \R_{\geq 0} (1, \Delta)  \subset \R \oplus 
M_{\R}$, where  $\deg (k, m) := k$ $(\forall k \in \Z, m \in M)$.  
Consider the graded 
semigroup ring 
\[ S_\Delta := \C[ M_\Delta] = \bigoplus_{k \geq 0} S^{(k)}_\Delta
 \subset \C[t_0, t_1^{\pm 1}, \ldots, t_d^{\pm 1}] , \]
 whose homogeneous components are $\C$-spaces of dimension 
 $|k\Delta \cap M|$. Note that $S^{(1)}_\Delta \subset 
 t_0 \C[t_1^{\pm 1}, \ldots, t_d^{\pm 1}]$ can be identified 
 with the  $\C$-vector 
 space $L(\Delta) \subset \C[t_1^{\pm 1}, \ldots, t_d^{\pm 1}] $ 
 spanned by all Laurent polynomials with the Newton 
 polytope $\Delta$.   

\begin{dfn}
A Laurent polynomial 
$f \in L(\Delta) \subset \C[t_1^{\pm 1}, \ldots, t_d^{\pm 1}]$  is called  {\em non-degenerate} if 
\[ t_0f, t_0t_1 \frac{\partial f}{\partial t_1},  \ldots,  t_0t_d \frac{\partial f}{\partial t_d} \in t_0L(\Delta) = S^{(1)}_\Delta \]
form a regular sequence  in 
$S_\Delta$. 
\end{dfn}

\begin{rem}
Let $A \subset \Delta \cap M$ be a finite subset which includes 
all vertices of the lattice polytope $\Delta$. Then non-degenerate 
Laurent polynomials 
\[ f({\bf t}) = \sum_{m \in A} a_m {\bf t}^m, \;\; a_m \in \C,  \]
form a Zariski dense open subset $U_A$ in 
the affine space $\C^{|A|}$ of 
all coefficients $\{ a_m \}_{m \in A}$. The subset $U_A$ 
can be explicitly 
defined by the non-vanishing of the principal $A$-determinant 
introduced by Gelfand, Kapranov and Zelevinski \cite{GKZ94}.
\end{rem} 

\begin{rem}
The Poincar{\'e} series of the graded ring $S_\Delta$  
\[ P(\Delta, t):= \sum_{k \geq 0 } |k \Delta \cap M|t^k \]
is a rational function of the form 
\[  P(\Delta, t) = \frac{\psi_0(\Delta) + \psi_1(\Delta) t + \cdots + \psi_d(\Delta)t^d}{(1-t)^{d+1} }, \]
where $\psi_0 (\Delta) = 0$ and $\sum_{i=0}^d \psi_i(\Delta) ={\rm Vol}_d (\Delta) = d!{\rm vol}(\Delta) \in \Z$. The coefficient 
$\psi_i(\Delta)$ is equal to the dimension of the $i$-th homogeneous component of the artinian ring
\[ S_\Delta/(t_0f,  t_0t_1 \frac{\partial f}{\partial t_1},  \ldots,  t_0t_d \frac{\partial f}{\partial t_d} ). \] 
\end{rem}

\begin{theorem} \cite[Remark 4.6]{DKh86} \label{dhE}
Let $f \in L(\Delta)$ be a non-degenerate Laurent polynomial and 
let $Z_{f, \Delta} \subset \T$ be the non-degenerate affine hypersurface defined by the equation $f =0$. Then 
the $E$-polynomial $E(Z_{f,\Delta}; u,1)$ has 
the following form: 
\[ E(Z_{f, \Delta}; u,1) = \frac{(u-1)^d - (-1)^d}{u} 
+ (-1)^{d-1} \sum_{i=1}^d \psi_i(\Delta) u^{i-1}. \]
\end{theorem}

\begin{rem} Explicit formulas for all Hodge-Deligne numbers
 $e^{p,q}(Z_f)$ as coefficients in  $E$-polynomial $E(Z_f; u,v)$  were obtained  by 
Danilov and Khovanski\v{i} only for simple 
$d$-dimensional lattice polytopes $\Delta$ \cite[Theorem 5.6]{DKh86}. The computations for arbitrary 
$d$-dimensional lattice polytopes $\Delta$  
are algorithmic  more complicated (see \cite{BB96,Sa20}). 
\end{rem}

In this paper we need explicit formulas for  
$E$-polynomials only in case when  $\Delta \subset M_\R$ is
a $d$-dimensional lattice simplex 
with vertices $v_0, v_1, \ldots, v_d$. We need also 
formulas for $E$-polynomials of affine hypersurfaces 
$Z_{f,\Delta_J}$ corresponding to faces $D_J$ 
of $\Delta$ parametrized by nonempty subsets $J \subseteq I:=\{0,1, \ldots, d \}$. If $J = \{j_1, \ldots, j_k \} \subset I$ is  a subset,  
 we denote $\Delta_J := {\rm conv}\{ v_{j_1}, \ldots, v_{j_k} \}$. 
 It is convenient to associate with 
a simplex $\Delta$  a finite abelian group $G(\Delta)$ 
defined as follows:

\begin{dfn}
Let $\Delta$ be a $d$-dimensional lattice simplex with 
vertices $v_0, v_1, \ldots, v_d$. We set $\widetilde{v_i}:= (1, v_i) \in \Z \oplus M$ $(0 \leq i \leq d)$ and 
denote by $M' \subset \Z \oplus M$ the sublattice of rank $d +1$
spanned  by linearly independent lattice vectors 
$\{ \widetilde{v_i} \}_{i \in I}$. Then 
\[ G(\Delta):= (\Z \oplus M)/M' \]
is a finite abelian group of order ${\rm Vol}_d(\Delta) = d!{\rm vol}(\Delta)$ and   
one can write an element $\widetilde{m} = (k, m) 
\in (\Z \oplus M)/M'$ as a unique $\Q$-linear combination 
\[ \widetilde{m} = (k, m)  = \sum_{i=0}^d \widetilde{\theta}_i(m)  \widetilde{v}_i , \]
where 
\[ (\widetilde{\theta}_0(m), \widetilde{\theta}_1(m), \ldots, \widetilde{\theta}_d(m)) \in [0,1)^{d+1}. \]
We set 
\[ {\rm age}(\widetilde{m}) = \sum_{i=0}^d   \widetilde{\theta}_i(m) = k , 
\;\; {\rm size}(\widetilde{m}) := {\rm age}(\widetilde{m}) + 
{\rm age}(-\widetilde{m}). \]
For any subset $J \subseteq I$ we define  
\[ G(\Delta_J):= \{ \widetilde{m} =(\widetilde{\theta}_0(m), \widetilde{\theta}_1(m), \ldots, \widetilde{\theta}_d(m))\; :\; 
 \widetilde{\theta}_j(m) = 0 \; \forall j \not\in J \}. \]
\end{dfn}

\begin{theorem} \label{E-simplex}
Let $\Delta \subset M_\R$ be a $d$-dimensional lattice 
simplex with $d+1$ lattice vertices $v_i$ $(i \in I)$, and 
let $J \subseteq I$ be a subset. 
Then 
\[ E(Z_{f,\Delta_J}; u,v) =   \frac{(uv-1)^{|J|-1} - 
(-1)^{|J|-1}}{uv} +  
\frac{(-1)^{|J|-2}}{uv} \sum_{0 \neq \widetilde{m} \in G(\Delta_J)}  v^{{\rm size}(\widetilde{m})} \left( \frac{u}{v} \right)^{{\rm age}(\widetilde{m})}. \]
\end{theorem}

\proof
First of all we note that 
\[ \psi_i(\Delta) = | \{ \widetilde{m} \in G(\Delta) \; :\; 
{\rm age}(\widetilde{m})  = i \}| \;\; \forall i \in \{ 0, 1, \ldots, d \},  \]
since one can choose a maximal 
regular sequence $s_0, s_1, \ldots, s_d \in S_\Delta^{(1)}$ corresponding to  $d+1$ monomials $\widetilde{v_0}, \widetilde{v_1}, 
\ldots,  \widetilde{v_d} \in M_\Delta$  and obtain a monomial 
$\C$-basis of $S_\Delta/(s_0,s_1, \ldots, s_d)$ that corresponds
to elements of $G(\Delta)$. Therefore, by \ref{dhE}, we obtain 
\[ E(Z_{f,\Delta}; u,1) =   \frac{(u-1)^{d} - 
(-1)^{d}}{u} +  
\frac{(-1)^{d-1}}{u} \sum_{0 \neq \widetilde{m} \in G(\Delta)}  
u^{{\rm age}(\widetilde{m})}. \]
Now we use the algorithm of Danilov and Khovanski\v{i} 
for simple polytopes $\Delta$ 
(see \cite[Sections 4 and 5]{DKh86}).   
We apply induction on dimension $d$ 
and assume that the formula in Theorem \ref{E-simplex} 
holds  for all proper faces $\Delta_J$  
of $\Delta$. Consider the Zariski 
closure  $\overline{Z}_{f,\Delta}$ of the 
affine hypersurface $Z_{f,\Delta}$  in the projective toric 
variety $\P_\Delta$. Then $\overline{Z}_{f,\Delta}$
 is quasi-smooth and the Hodge structure 
in $H^*_c(\overline{Z}_{f,\Delta})$ is pure, i.e., only 
elements $\widetilde{m} \in G(\Delta)$ with the maximal ${\rm size}(\widetilde{m})
= d+1$ can contribute to $H_c^{d-1}(\overline{Z}_{f,\Delta})$. 
It remains to apply the additivity and the induction hypothesis to   
the equality 
\[ E(\overline{Z}_{f,\Delta}; u,v) = 
\sum_{J \subset I}  E(Z_{f,\Delta_J}; u,v). \]
\endproof

\begin{coro} \label{E-cyclic}
Let $\overline{w} = (w_0, w_1, \ldots, w_d) \in \Z^{d+1}_{>0}$ be 
a well-formed weight vector and lel $\Delta \subset M_\R$ be 
the $d$-dimensional lattice 
simplex whose vertices $v_0, v_1, \ldots, v_d$ generate 
$M$ and satisfy the relation $\sum_{i=0}^d w_i v_i = 0$. 
For any subset 
$J \subset I = \{0,1, \ldots, d \}$ consider 
the subgroup 
\[ G_J := \{ l \in  \Z/w\Z \; :\;  \widetilde{\theta}_j(l) =0 \;\; \forall 
j \not\in J \} \subseteq G_I= G_0 =  \Z/w\Z.  \]
Then 
\[ E(Z_{f,\Delta_J}; u,v) =   \frac{(uv-1)^{|J|-1} - 
(-1)^{|J|-1}}{uv} +  
\frac{(-1)^{|J|-2}}{uv} \sum_{0 \neq l \in G_J}  v^{{\rm size}(l)} \left( \frac{u}{v} \right)^{{\rm age}(l)}. \]
\end{coro}

\proof
The statement follows immediately from Theorem \ref{E-simplex}.  
Note that this particular case was considered by Corti and Golyshev in 
\cite[Theorem 1.3]{CG11}. 
\endproof

\section{Stringy $E$-functions of Calabi-Yau hypersurfaces}

Let $M \cong \Z^d$ be a lattice of rank $d$, $N:= {\rm Hom}(M, \Z)$ the dual lattice, and 
$\langle *, * \rangle \, :\, M \times N \to \Z$ the 
natural pairing. We set  
$M_\R := M \otimes \R$ and  $N_\R:= N \otimes \R$. 
By a  {\em lattice polytope} $\Delta \subset M_\R$ we mean
a convex hull of 
finitely many lattice points $A \subset M$, i.e., 
$\Delta = {\rm conv}(A)$. For a $d$-dimensional 
lattice polytope $\Delta$ we consider the piecewise linear function 
\[ {\rm ord}_\Delta\, : \, N_\R \to \R, \;\;  {\rm ord}_\Delta(y):= 
\min_{x \in \Delta} \langle x, y \rangle. \]
The domains of linearity of ${\rm ord}_\Delta$  form a complete 
rational polyhedral fan $\Sigma_\Delta$ in $N_\R$ which is called 
the {\em normal  
fan of $\Delta$} \cite{CLS11}.  The normal fan $\Sigma_\Delta$  consists of cones $\sigma_\Theta$ 
which $1$-to-$1$ correspond to faces $\Theta \preceq \Delta$: 
\[ \sigma_{\Theta}:= \{ y \in N_\R\; :  \;   
\langle x, y \rangle  = {\rm ord}_\Delta(y) \;\; \forall x \in \Theta\}. \]
 
\begin{dfn}
Let $\Delta \subset M_\R$ be a $d$-dimensional 
lattice polytope. The subset 
\[    F(\Delta) := 
\{ x \in M_\R \, : \, \langle x, n \rangle \geq 
{\rm ord}_\Delta(n) + 1 \; \; \forall n \in N \setminus \{0 \} \} \subset \Delta. \]
is called the {\em Fine interior} of $\Delta$.   
\end{dfn}

\begin{rem}
The Fine interior of a lattice polytope $\Delta$ is either empty, or
a rational polytope of dimension $\leq d$. 
The convex hull ${\rm conv}(\Delta^\circ \cap M)$ is always contained in $F(\Delta)$. If $d =2$, then $F(\Delta) = 
{\rm conv}(\Delta^\circ \cap M)$ \cite[Prop.2.9]{Bat17}. There exist exactly 
$9$ lattice  polytopes $\Delta$ of dimension $d =3$ such that
$F(\Delta) \neq \emptyset$, but $ \Delta^\circ \cap M = \emptyset$ \cite[Appendix B]{BKS19}. 
\end{rem}

\begin{theorem} \cite[Theorem 2.23]{Bat17}
A  non-degenerate affine hypersurface $Z_\Delta \subset \T
:= {\rm Hom}(M, \C^*)$  
defined by a Laurent polynomial $f \in \C[M]$ 
with the Newton polytope $\Delta$ admits a Calabi-Yau 
compactification $X_\Delta$ if and only if the Fine interior 
of $\Delta$ is a single lattice point. 
\end{theorem}

\begin{dfn}
A $d$-dimensional lattice polytope $\Delta \subset M_\R$ 
containing $0 \in M$ in its interior is called {\em canonical 
Fano polytope} if $0 = \Delta^\circ \cap M$, i.e., $0$  is the unique interior lattice point in 
$\Delta$. 
\end{dfn} 
 
\begin{rem} If $F(\Delta) = \{0 \}$, then $\Delta$ is a canonical 
Fano polytope, but the converse is not true if $\dim \Delta \geq 3$. 
All $3$-dimensional canonical Fano polytopes are classified 
by Kasprzyk \cite{Kas10}. Among them 
 exist exactly $9089$ canonical 
Fano polytopes $\Delta$ such that $\dim F(\Delta) \geq 1$ \cite{BKS19}. 
\end{rem}

It is elementary to show the following:

\begin{prop}
Let $\Delta$ be a $d$-dimensional canonical Fano polytope. 
Denote by $\Delta^* \subset N_\R$ the dual rational polytope
\[ \Delta^*:= \{ y \in N_\R \; : \; \langle x, y \rangle \geq -1 \; \; \forall x \in \Delta\}. \] 
Then $F(\Delta) = \{ 0 \}$ if and only if 
\[ [ \Delta^*] := {\rm conv}(\Delta^* \cap N) \subset N_\R \] 
is also a $d$-dimensional canonical Fano polytope. 
\end{prop}

\begin{coro}
Let $\overline{w} \in \Z^{d+1}_{>0}$ be a weight vector. 
A general hypersurface $X_w \subset \P(\overline{w})$ of 
degree $w = \sum_{i=0}^d w_i$ is a Calabi-Yau variety
if and only if 
$\overline{w}$ has $IP$-property.
\end{coro}

\begin{coro} \label{CY-affine}
Let $\overline{w} \in \Z^{d+1}_{>0}$ be a weight vector. Consider 
$d$-dimensional simplex $\Delta_{\widetilde{w}}^* 
= {\rm conv}(v_0, v_1, \ldots,v_d) \subset M_\R$ 
such that $M = \sum_i \Z v_i$ and $\sum_i w_i v_i =0$. 
Then $F(\Delta_{\widetilde{w}}^*) = \{ 0 \}$ if and only 
if $\overline{w}$ has $IP$-property.
\end{coro}

\begin{rem}
Assume that $\Delta$ and $[\Delta^*]$ are canonical Fano polytopes, i.e., $F(\Delta) = \{0 \}$. There exist many ways to obtain a Calabi-Yau compactifications of a non-degenerate affine hypersurface 
$Z_\Delta\subset \T$. For constructing a Calabi-Yau 
compactification  one has to choose 
a finite subset $B \subset [\Delta^*] \cap N$ such that 
$\nabla_B:= {\rm conv}(B)$ is a canonical Fano polytope 
and consider the fan $\Sigma^B$ in $N_\R$ with cones $\R_{\geq 0} \nabla$ spanned by  
proper faces $\nabla \prec \nabla_B$ and   
$\{0 \} \in N_\R$. Then a Calabi-Yau compactification $
X_\Delta^B$ of $Z_\Delta$ is the Zariski 
closure of $Z_\Delta$ in 
$\Q$-Gorenstein Fano toric variety $\P_B$ associated with the 
fan $\Sigma^B$ \cite[Theorem 1]{ACG16}.  
\end{rem} 

\begin{theorem} \cite[Theorem 4.10]{Bat17} \label{E-str}
Let  $\Delta \subset M_\R$ 
be a $d$-dimensional lattice polytope with 
$F(\Delta) = \{ 0 \}$. Then 
the stringy $E$-function of a Calabi-Yau compactification 
$X_\Delta$ of a non-degenerate affine 
hypersurface $Z_\Delta \subset \T$ with the Newton 
polytope $\Delta$ can be computed by the formula
\begin{align*} 
E_{\rm str}(X_\Delta; u,v)= \sum_{ \Theta \preceq \Delta \atop 
\dim \Theta \geq 1 }  E(Z_\Theta; u,v ) 
(uv - 1 )^{d - \dim \Theta} \sum_{ n \in \sigma_{\Theta}^\circ} (uv)^{{\rm ord}_\Delta(n)}, 
\end{align*}
where $E(Z_\Theta; u,v )$ denotes the 
 $E$-polynomial 
of the affine hypersurface $Z_\Theta \subset \T_\Theta$
corresponding to a face $\Theta \preceq \Delta$, and $\sigma_\Theta^\circ$ is the relative interior of the $(d -\dim \Theta)$-dimensional cone $\sigma_\Theta \in \Sigma_\Delta$. 
\end{theorem}

\begin{rem}
We note that the combinatorial  formula in Theorem \ref{E-str} 
depends only on $\Delta$ and it does not depend on the choice 
of a Calabi-Yau compactification $X_\Delta^B$. 
\end{rem}

We specialize the formula in Theorem \ref{E-str} to the case 
when $$\Delta = \Delta_{\widetilde{w}}^* 
= {\rm conv}(v_0, v_1, \ldots,v_d) \subset M_\R$$ for a  
weight vector $\overline{w}$ with $IP$-property. 
In this case 
 we need the following standard combinatorial fact about the set $L^{\overline{w}}(k)$ of nonnegative 
integral solutions $\overline{u} = (u_0, u_1, \ldots, u_d) \subset \Z^{d+1}_{\geq 0}$ to the linear diophantine equation:
\[ \sum_{i=0}^d w_i u_i  = \langle \overline{w},  \overline{u} \rangle = kw,  \; \; k\in \Z_{>0}.  \]

\begin{prop} \label{sum}
For any subset $J \subseteq I = \{0,1, \ldots, d\}$ denote by $L_J^{\overline{w}}(k)$ the subset in $L^{\overline{w}}(k)$:
\[ L_J^{\overline{w}}(k) := \{ \overline{u} \in L^{\overline{w}}(k) \, : \, 
u_j = 0 \Leftrightarrow j \in J \}. \]
Then 
\[ \sum_{k >0} |L_J^{\overline{w}}(k) | t^{-k} = \left[ \prod_{j \not\in J} \frac{1}{t^{q_j} -1} \right]_{int}. \]
\end{prop}
\proof
Consider $d+1$ variables 
$t_0, t_1, \ldots, t_d$ and associate with a
lattice point $\overline{u} \in \Z^{d+1}_{\geq 0}$ 
the monomial ${\bf t}^{\overline{u}} =\prod_i t_i^{u_i}$. 
A lattice point $\overline{u} \in \Z^{d+1}_{\geq 0}$ belongs 
to $L^{\overline{w}}(k)$ if and only if by setting $t_i = t^{-q_i}$ 
$( 0 \leq i \leq d)$ we 
get a monomial in $t$ with the integral exponent $-k$, i. e.,  
\[ \prod_{i=0}^d (t^{-q_i})^{u_i} = t^{-k}.  \]
Therefore, 
\[ \sum_{k >0} |L_J^{\overline{w}}(k) | t^{-k} = 
\left[ \prod_{j \not\in J} \left( t^{-q_j} + t^{-2q_j} + \cdots \right) \right]_{int} = 
\left[ \prod_{j \not\in J} \frac{1}{t^{q_j} -1} \right]_{int}. \]
\endproof

\begin{coro} \label{sum-sigma}
Let $\Delta = \Delta_{\overline{w}}^* = {\rm conv}(v_0, v_1, \ldots, v_d)$ be $d$-dimensional simplex  whose vertices a primitive vectors satisfying  
the relation $\sum_{i=0}^d w_i v_i =0$. Then for any $J \subset I$
one has 
\[ \sum_{ n \in \sigma_J^\circ} (uv)^{{\ord}_\Delta(n)} =
\left[ \prod_{j \not\in J} \frac{1}{(uv)^{q_j} -1} \right]_{int}. \]  
\end{coro}

\proof  The function ${\rm ord}_\Delta$ is linear 
on the $(d+1 - |J|)$-dimensional simplicial cone $\sigma_J^\circ$ and 
it 
has value $(-1)$ on every vertex of the rational dual simplex 
$\Delta^* = \Delta_{\overline{w}}$. 
The set of lattice points $n \in \sigma_J^\circ$ 
such that ${\rm ord}_\Delta (n) = -k$
can be identified with the subset $L_J^{\overline{w}}(k)$ 
of nonnegative integral solutions of the linear 
diophantine equation $\sum_i w_i u_i = w$.  Now the 
statement follows from  \ref{sum}.
\endproof

\begin{theorem} \label{E-str-final}
Let $\overline{w} \in \Z^{d+1}$ be a weight vector with 
$IP$-property.  Then a Calabi-Yau compactification 
$X_{\overline{w}}^\vee$ of a non-degenerate hypersurface
$Z_{\overline{w}} \subset \T_{\overline{w}}^*$ defined by a Laurent polynomial with Newton polytope $\Delta^*_{\overline{w}}$ has
the following stringy $E$-function 
\[ E_{\rm str}(X_{\overline{w}}^\vee; u,v) = \sum_{J \subseteq I \atop |J| \geq 2} E(Z_{\overline{w},J}; u,v) (uv-1)^{d+1 - |J|} 
\left[ \prod_{j \not\in J} \frac{1}{(uv)^{q_j} -1} \right]_{int}, \]
where 
\[ E(Z_{\overline{w},J}; u,v) =    \frac{(uv-1)^{|J|-1} - 
(-1)^{|J|-1}}{uv} +  
\frac{(-1)^{|J|-2}}{uv} \sum_{0 \neq l \in G_J}  v^{{\rm size}(l)} \left( \frac{u}{v} \right)^{{\rm age}(l)}. \] 
\end{theorem}

\proof The statement follows immediately from Theorem \ref{E-str}
using \ref{sum-sigma} and \ref{E-cyclic}. 
\endproof

In Example \ref{E-str} below we illustrate the combinatorial 
formula from Theorem \ref{E-str-final}. 

\begin{exam}
Take the weight vector $\overline{w} = (1,5,12,18)$. It has 
$IP$-property. The 
affine hypersurface  $Z_f \subset (\C^*)^3$ defined by the equation 
\[ f_{\overline{w}}^0({\bf t}) = \frac{1}{t_1^5 t_2^{12} t_3^{18}} +
t_1 + t_2 + t_3 = 0 \]
is non-degenerate.   
The Newton polytope of $f_{\overline{w}}^0({\bf t})$ is 
a $3$-dimensional lattice simplex $\Delta$ with $F(\Delta) = \{ 0 \}$
and ${\rm Vol}_3(\Delta)= 36 = w$. By \ref{CY-affine}, 
$Z_f$ admits a
Calabi-Yau compactification $X_{\Delta}$. We compute 
its stringy $E$-function $E_{\rm str}(X_\Delta; u,v)$ using 
Theorem \ref{E-str}. 

Note that the sum over faces $\Delta_J\preceq \Delta$ $(\dim \Delta_J \geq 1)$ 
in the formula \ref{E-str-final} consists of $1 + 4 + 6 = 11$ terms ($\dim 
\Delta_J \in \{ 3,2,1 \}$). By \ref{E-cyclic}, we obtain 
\[ E(Z_f; u,v)= \left( (uv)^2 - 3uv + 3 \right) + (u^2 +10 uv + v^2) + 3(u+v) + 6(u+v) + 5 = \]
\[ =  (uv)^2 + u^2 + 7uv + v^2 + 9(u+v) +8.  \] 
The terms 
corresponding to four $2$-dimensional faces $\Delta_J \prec 
\Delta$ are
\[ (uv-w), \; (uv-2) \frac{uv-1}{(uv)^5-1}, \; (uv - 3(u +v ) -7), 
\;  (uv - 6(u+v) -7). \] 
The terms corresponding to six $1$-dimensional faces 
$\Delta_J \prec \Delta$ are 
\[ \frac{(1 + 7uv + 7(uv)^2 + 
7(uv)^3 + 7(uv)^4 +7(uv)^5)(uv-1) }{(uv)^5 -1}, \;  (1 + 2uv), \;\;
(1 + uv),  \]
\[ 6, \; \; \frac{(1+ (uv)^3)(uv -1)}{(uv)^5 -1}, \; \; 
\frac{(1+ (uv)^2 + (uv)^4)(uv -1)}{(uv)^5 -1}. 
\]

%$E(\Delta_{(1)}) \times $: $(uv - 2) \times 1 = uv -2  $

%$E(\Delta_{(5)}) \times $: $(uv - 2)\times \frac{uv-1}{(uv)^5-1}$

%$E(\Delta_{(12)}) = (uv - 3(u+v) -7) \times 1  $ 

%$E(\Delta_{(18)}) = (uv - 6(u+v) -7) \times 1 $. 

%$E(\Delta_{(1),(5)})   = 1 \times \frac{(1 + 7uv + 7(uv)^2 + 
%7(uv)^3 + 7(uv)^4 +7(uv)^5)(uv-1) }{(uv)^5 -1}$

%$E(\Delta_{(1),(12)}) = 1 \times (1 + 2uv)$ 

%$E(\Delta_{(1),(18)}) = 1 \times (1 + uv)$  

%$E(\Delta_{(12),(18)})  = 6 \times 1 $

%$E(\Delta_{(5),(18)} = 1 \times \frac{(1+ (uv)^3)(uv -1)}{(uv)^5 -1}$ 

%$E(\Delta_{(5),(12)})  = 1 \times  \frac{(1+ (uv)^2 + (uv)^4)(uv -1)}{(uv)^5 -1}$ 
The stringy Euler number $E_{\rm str}(X_\Delta; 1,1)$ equals
\[  36 - \frac{1}{5} - 1 - 12 - 18 + 1 \cdot \frac{36}{5} + 1 \cdot 2 + 1 \cdot 3  + 6 \cdot 1 + 1 \cdot \frac{2}{5} + 1 \cdot \frac{3}{5} =24. \]
Note that the sum of all four terms 
with the denominator $(uv)^5 -1$
 equals $1 + 7uv$. Thus, we obtain
%$\Delta_{(1),(5)} + \Delta_{(5)} + \Delta_{(5),(12)}  + 
%\Delta_{(5), (18)}= 1 + 7uv$. 
\[ E_{\rm str}(X_\Delta; u,v) 
= (uv)^2 + u^2 + 7uv + v^2 + 9(u+v) +8  + (7uv +1)  + \]
\[ + 
 (uv -2) + (uv - 3(u+v)  -7) + (uv - 6(u +v) - 7) + 6 + (2 + 3 uv) = 
 \]
 \[ = (uv)^2 + u^2 + 20uv + v^2 + 1.   \]  
 This is $E$-polynomial of a $K3$-surface. 
  \end{exam}

\section{The formula of Vafa}

Let $W \in R:= \C[z_0,z_1, \ldots, z_d]$ be a non-degenerate 
weighted homogeneous polynomial with 
weights $w_i:= \deg z_i$ $(0 \leq i \leq d)$ and
 $\deg W = w =\sum_{i=0}^d w_i$ , i.e., 
\[ W(\lambda^{w_0}z_0, \ldots, \lambda^{w_d}z_d) = \lambda^w W(z_0, \ldots, z_d) \;\;\forall \lambda \in \C^* \]
and  the common zero of all partial derivatives 
\[ W_i' :=\frac{\partial W}{\partial z_i} , \; \; (0 \leq i \leq d) \]
is only the origin $0 \in \C^{n+1}$. Then $d+1$ homogeneous 
polynomials  $W_i$ of degree $w-w_i$ $(0 \leq i \leq d)$ form a regular sequence in the graded ring $R$ and the quotient $\overline{R}:= R/\langle W_0', \ldots, W_d' \rangle$ is a graded artinian ring with   
homogeneous components $\overline{R}_m$ 
 such that 
\[ P(W, t) = \sum_{m \geq 0} \dim_\C \overline{R}_m t^m  = \prod_{j =0}^d \left( \frac{1 - t^{w -w_i}}{1 - t^{w_i}} \right)  \]
and 
\[ \mu:= \dim_\C \overline{R} = \prod_{j=0}^d  \left( \frac{{w -w_i}}{{w_i}} \right)    \]
is the {\em Milnor number} of $W$. 
As above, we set 
\[  (q_0, q_1, \ldots, q_d):= \frac{1}{w}(w_0, w_1, \ldots, w_d) \]
and consider the cyclic group $G_0 = \Z/w\Z$ 
that acts on $\C^{n+1}$ by diagonal matrices 
\[ \exp (2\pi  i \theta(l)) = {\rm diag}(e^{2\pi i \theta_0(l)}, \ldots, e^{2\pi i \theta_d(l)}  ), \;\;  l \in \Z/w\Z.  \]
For any element $l \in \Z/l\Z$ we define the 
subspace 
\[ U^{(l)}  := \{ (z_0, z_1, \ldots, z_n) \in \C^{n+1} \, :  
\, z_j =0 \;  \, \forall  
lq_j \in \Z\}. \]
Then the restriction $W^{(l)}$ 
of the polynomial $W$ to the subspace 
$U^{(l)} \subset \C^{d+1}$ is also non-degenerate and 
therefore we obtain another 
polynomial 
\[ P(W^{(l)}, t):= \prod_{lq_j \in \Z} \left( \frac{1 - t^{w -w_j}}{1 - t^{w_j}} \right).   \]

An explicit method for computing the orbifold  Hodge
numbers $h^{p,q}_{\rm ord}(X_w)$ of quasi-smooth Calabi-Yau hypersurfaces $X_w$
 was first suggested by 
Vafa \cite{Vaf89}. He gave an explicit 
formula for the orbifold Euler 
\begin{align} \label{vafa0}
\chi_{\rm orb}(X_w) = \frac{1}{w} 
\sum_{l,r=0}^{w-1} \, \, \prod_{0 \leq i \leq d \atop lq_i,rq_i \in \Z}
\left(1 - \frac{1}{q_i} \right).
\end{align}
 The method  
for computing the orbifold Hodge numbers of $X_w$ 
was illustrated by Vafa in some examples \cite[p.1182-1183]{Vaf89}, but an exact combinatorial mathematical formula based on his ideas appeared 
later 
in the paper of Kreuzer, Schimmrigk and Skarke \cite[Formula (6)]{KSS92}(see also Klemm, Lian, Roan and Yau \cite[Formula (3.2)]{KLRY98}):

\begin{align} \label{vafa1}
P(t, \overline{t}) = \sum_{0 \leq l < w} \left[ 
\prod_{\widetilde{\theta}_i(l) = 0} 
\frac{1-(t\overline{t})^{1 -q_i}}{1-(t\overline{t})^{q_i}} \prod_{\widetilde{\theta}_i(l) \neq 0} 
(t \overline{t})^{\frac{1}{2} - q_i}\left(  \frac{t}{\overline{t}}\right)^{\widetilde{\theta}_i(l) - \frac{1}{2}} \right]_{int},  
\end{align}
where 
\[ \theta(l) = l(q_0, \ldots, q_d), \;\; \widetilde{\theta}(l)= 
(\widetilde{\theta}_0(l), \ldots, \widetilde{\theta}_d(l)) \in [0,1)^{d+1}, \] 

\begin{rem} \label{mirr-pq}
It is important to note that the polynomial $P(t, \overline{t}) = \sum_{p,q} h^{p,q} t^p \overline{t}^q $ in the formula
(\ref{vafa1})
has nonnegative integral coefficients $h^{p,q}$ which are 
not the orbifold Hodge numbers $h^{p,q}_{\rm orb}(X_w)$ of the
quasi-smooth Calabi-Yau hypersurface, but 
the Hodge numbers $h^{d-1-p,q}$ of its mirror. This can be seen
more explictly in examples below. 
In particular, for Calabi-Yau $3$-folds $X$ one obtains  
\[ P(t, \overline{t}) = (1 + t^3)(1 + \overline{t}^3) + 
 h^{2,1}(X) (t \overline{t} + (t\overline{t})^2) 
+ h^{1,1}(X) (t (\overline{t})^2 + t^2 \overline{t}).  \]
\end{rem}

\begin{exam}
Let $d = 4$, $w = 5$ and $\overline{w} := (1,1,1,1,1)$. 
Then for $l =0$ we have 
\[ P_0(t, \overline{t}) =
\left[ \prod_{j =0}^4 \left( \frac{1-(t \overline{t})^{4/5}}
{1-(t\overline{t})^{1/5}} \right) \right]_{int}  
=\left[ (1 + (t \overline{t})^{1/5} + (t \overline{t})^{2/5} + 
(t \overline{t})^{3/5})^5 \right]_{int} = \]
\[ 1 + 101 t\overline{t} + 101 (t\overline{t})^2 +  (t\overline{t})^3. \] 
For the remaining values $l \in \{1,2,3,4 \}$,  
one obtains four monomials 
\[ t^3, t^2 \overline{t},  t \overline{t}^2, \overline{t}^3. \]
Since a Calabi-Yau quintic $3$-fold $X_5 \subset P^4$ 
has Hodge numbers $h^{1,1}(X_5) =1$ and $h^{2,1}(X_5) 
= 101$,  the formula (\ref{vafa1}) produces the Hodge numbers 
of the mirror of $X_5$.  
\end{exam}
 
We need a modified version of the 
formula (\ref{vafa1}): 

\begin{prop} Let 
\[ Q(u,v) :=  \left[ \sum_{0 \leq l < w} \left( \prod_{\theta_i(l) \in \Z} \frac{1 - (uv)^{1 -q_i}}{1-(uv)^{q_i}} \right)  \cdot \prod_{\theta_i(l) \not\in \Z} (uv)^{\frac{1}{2} - q_i} \left( \frac{u}{v} \right)^{\widetilde{\theta}_i(l) -\frac{1}{2}} \right]_{int}. \] 
Then one has
\[ Q(u,v) =   \frac{1}{uv}  \sum_{0 \leq l < w}  
 \left[ \prod_{\widetilde{\theta}_i(l) \in \Z} \frac{(uv)^{q_i}  - uv}{1 - (uv)^{q_i}}  \right]_{int}
\cdot    v^{{\rm size}(l)} 
\left( \frac{u}{v} \right)^{{\rm age}(l)}.   \]
\end{prop}

\proof 
By  $ \sum_{i=0}^d q_i =1$,  we have
\[ \prod_{\widetilde{\theta}_i(l) \not\in \Z} (uv)^{- q_i}  = 
\frac{1}{uv} \prod_{\widetilde{\theta}_i(l) \in \Z} (uv)^{ q_i}. \] 
Using 
\[ \prod_{\widetilde{\theta}_i(l) \not\in \Z} (uv)^{1/2 - q_i}  = 
(uv)^{\frac{{\rm size}(l)}{2} } \prod_{\widetilde{\theta}_i(l) \not\in \Z} (uv)^{ - q_i}, \]
we obtain 
\[ Q(u,v) =  \left[ \sum_{0 \leq l < w} \left( \prod_{\theta_i(l) \in \Z} \frac{1 - (uv)^{1 -q_i}}{1-(uv)^{q_i}} \right)  \cdot \prod_{\theta_i(l) \not\in \Z} (uv)^{\frac{1}{2} - q_i} \left( \frac{u}{v} \right)^{\widetilde{\theta}_i(l) -\frac{1}{2}} \right]_{int}   = \] 
 \[  =\sum_{0 \leq l \leq w-1} \frac{1}{uv}  \left[   \left( \prod_{\widetilde{\theta}_i(l) \in \Z} 
\frac{(uv)^{q_i}  - uv}{1 - (uv)^{q_i}} \right) (uv)^{\frac{{\rm size}(l)}{2}} 
 \left( \frac{u}{v} 
\right)^{- \frac{{\rm size}(l)}{2}} 
\cdot \prod_{\widetilde{\theta}_i(l) \not\in \Z} \left( \frac{u}{v} 
\right)^{\widetilde{\theta}_i(l) } \right]_{int}  = \]

\[ \frac{1}{uv} \sum_{0 \leq l < w}  
\left[  \prod_{\widetilde{\theta}_i(l) \in \Z} \frac{(uv)^{q_i}  - uv}{1 - (uv)^{q_i}} \right]_{int}   \cdot 
\left(  v^{{\rm size}(l)} \prod_{\widetilde{\theta}_i(l) \not\in \Z} 
\left( \frac{u}{v} \right)^{\widetilde{\theta}_i(l)} \right) =  \]

\[ =\frac{1}{uv}  \sum_{0 \leq l < w}  
 \left[ \prod_{\widetilde{\theta}_i(l) \in \Z} \frac{(uv)^{q_i}  - uv}{1 - (uv)^{q_i}}  \right]_{int}
\cdot    v^{{\rm size}(l)} 
\left( \frac{u}{v} \right)^{{\rm age}(l)}.  \]
\endproof

\begin{coro}  \label{mirr-trans}
Assum that $\overline{w} \in \Z^{d+1}_{>0}$ is a transverse weight vector and $X_w \subset \P(\overline{w})$ is 
a quasi-smooth Calabi-Yau hypersurface. Then for the mirror 
stringy $E$-polynomial of $X_w$ we have
\[ (-u)^{d-1} 
E_{\rm orb}(X_w; u^{-1},v) =  \frac{1}{uv}  \sum_{0 \leq l < w}  
 \left[ \prod_{\widetilde{\theta}_i(l) \in \Z} \frac{(uv)^{q_i}  - uv}{1 - (uv)^{q_i}}  \right]_{int}
\cdot    (-v)^{{\rm size}(l)} 
\left( \frac{u}{v} \right)^{{\rm age}(l)}. \] 
\end{coro}

\begin{exam} Consider the transverse weight vector
$\overline{w}= (1,5,12,18)$. Then $w = \sum_i w_i = 36$. 
There are twelve  elements $l$ in $\Z/36\Z$ 
acting on $\C^4$ with size $4$:
one element of age $1$, ten elements of age $2$, and one element of age $3$. This gives polynomial 
$u^2 + 10 uv + v^2$. 
The element $l =0$ determines the term 
\[ \left[  \left( \frac{1-(uv)^{\frac{35}{36}}}{1- (uv)^{\frac{1}{36}}} \right) \cdot \left( \frac{1-(uv)^{\frac{31}{36}}}{1- (uv)^{\frac{5}{36}}} \right)  \left( \frac{1-(uv)^{\frac{2}{3}}}{1- (uv)^{\frac{1}{3}}} \right) 
\left( \frac{1-(uv)^{\frac{1}{2}}}{1- (uv)^{\frac{1}{2}}} \right)  \right]_{int} = 1 + 10 uv + (uv)^2. 
\]
There exist five elements $l \in \Z/36\Z$ 
of size $2$ and age $1$ with the corresponding 
term 
\[ \left[ (uv)^{\frac{5}{6}}\left( \frac{1-(uv)^{\frac{2}{3}}}{1- (uv)^{\frac{1}{3}}} \right) 
\left( \frac{1-(uv)^{\frac{1}{2}}}{1- (uv)^{\frac{1}{2}}} \right) 
\right]_{int} =    \left[ (uv)^{\frac{5}{6} }(1 + (uv)^{\frac{1}{3}} )\right]_{int} =0.   \] 
There exist twelve elements $l \in \Z/36\Z$ 
of size $3$ and age $1$ or $2$ with the corresponding 
term  
\[  \left[ (uv)^{\frac{1}{2}}
\left( \frac{1-(uv)^{\frac{1}{2}}}{1- (uv)^{\frac{1}{2}}} \right) 
\right]_{int}\cdot    (-v)
{u}^{{\rm age}(l)-1} = \left[ (uv)^{\frac{1}{2} } \cdot  1 \right]_{int}\cdot    (-v) {u}^{{\rm age}(l)-1} =0.   \] 
There exist six elements $l \in \Z/36\Z$ 
of size $3$ and age $1$ or $2$ with the corresponding 
term  
\[  \left[ (uv)^{\frac{1}{3}} \left( \frac{1-(uv)^{\frac{2}{3}}}{1- (uv)^{\frac{1}{3}}} \right) \right]_{int}\cdot    (-v) {u}^{{\rm age}(l)-1} =  \left[ (uv)^{\frac{1}{3}} 
(1 +  (uv)^{\frac{1}{3}} \right]_{int}\cdot   (-v) {u}^{{\rm age}(l)-1}   =0.   \] 
By summing all $36$ terms, we obtain 
\[ (-u)^2E_{\rm orb}(X_w; u^{-1},v) =  1 + 10 uv + (uv)^2 +
 u^2 + 10 uv + v^2 = 1 +  u^2 + 20 uv + v^2  + (uv)^2. \]
\end{exam}

\begin{exam}  Take the transverse weight vector $ \overline{w} = (1,1,2, 2,2)$. Then  $w = 8$ and we obtain the action of 
$\Z/8\Z$ on $\C^5$ with charges $(1/8, 1/8, 1/4, 1/4, 1/4)$. 
There exist $6$ elements $l \in \Z/8\Z$ of size $5$ (one of age $1$, two of age $2$, two of age $3$  and one of age $4$). This gives a polynomial $- v^3 - 2 v^2u - 2v u^2 - u^3$.  
The element $l =0$ in $\Z/8\Z$ determines 
the term 
\[ \left[  \left( \frac{1-(uv)^{\frac{7}{8}}}{1- (uv)^{\frac{1}{8}}} \right)^2 \cdot  \left( \frac{1-(uv)^{\frac{3}{4}}}{1- (uv)^{\frac{1}{4}}} \right)^3  \right]_{int} = 1 + 83 uv + 83 (uv)^2 + (uv)^3. 
\]
It remains to consider 
the element of order $2$ in $l \in \Z/8\Z$. It  has size $2$ and 
age $1$.  By putting $t = (uv)^{\frac{1}{4}}$ in the equality
\[ (1 + t + t^2)^3 = 1 + 3t+ 6 t^2 + 7 t^3 + 6 t^4 + 3 t^{5} + 
t^{6},  \]
we  obtain 
the corresponding term 
\[ \left[   (uv)^{\frac{3}{4}} \cdot  \left( \frac{1-(uv)^{\frac{3}{4}}}{1- (uv)^{\frac{1}{4}}} \right)^3  \right]_{int} = 3uv + 3 (uv)^2. 
\]
Therefore, 
we obtain 
\[   (-u)^3E_{\rm orb}(X_w; u^{-1},v) = 1  + 
86 uv - v^3 - 2 v^2u - 2v u^2 - u^3 + 86 (uv)^2  + (uv)^3. \] 
In particular, the Euler number of $X_8 \subset \P(1,1,2,2,2)$ 
equals  $-168$. 
\end{exam}

%\section{E-polynomials}

\section{The proof}

 \noindent
{\bf Proof of Theorem 2.1} 
In order to prove the required equality 
\[ E_{\rm str} (X_{\overline{w}}^\vee; u, v) = 
\frac{1}{uv}  \sum_{0 \leq l < w}  
 \left[ \prod_{\widetilde{\theta}_i(l) = 0} \frac{(uv)^{q_i}  - uv}{1 - (uv)^{q_i}}  \right]_{int}
\cdot   (- u)^{{\rm size}(l)} 
\left( \frac{v}{u} \right)^{{\rm age}(l)}   \]
we split its left hand side and right hand side
 into a sum over elements $l \in G_0 = \Z/w\Z$. We 
set
\[ P_{\overline{w}}^{(l)}(u,v):= \frac{1}{uv}\left[ \prod_{\widetilde{\theta}_i(l) = 0} \frac{(uv)^{q_i}  - uv}{1 - (uv)^{q_i}}  \right]_{int}
\cdot   (- u)^{{\rm size}(l)} 
\left( \frac{v}{u} \right)^{{\rm age}(l)}, \;\; l \in \Z/w\Z, \]
and use  the formula in Theorem \ref{E-str-final} to write  
\[ E_{\rm str} (X_{\overline{w}}^\vee; u, v)   = \sum_{0 \leq l < w} E^{(l)}_{\rm str}(X_{\overline{w}}^\vee; u, v), \;\; 
\]
where 
\[ E^{(0)}_{\rm str}(X_{\overline{w}}^\vee; u, v):= \sum_{J \subseteq I \atop |J| \geq 2}  \frac{(uv-1)^{|J|-1} - 
(-1)^{|J|-1}}{uv} \cdot (uv -1)^{d + 1 -|J|} \left[ \prod_{j \not\in J} \frac{1}{(uv)^{q_j} -1} \right]_{int} \] 
and 
\[ E_{\rm str}^{(l)}(X_{\overline{w}}^\vee; u,v): = \sum_{J \subseteq I \atop l \in G_J }  
\frac{ (-1)^{|J|}v^{{\rm size}(l)}}{uv} \left( \frac{u}{v} \right)^{{\rm age}(l)} \cdot (uv -1)^{d + 1 -|J|} \left[ \prod_{j \not\in J}  \frac{1}{(uv)^{q_j} -1}
 \right]_{int}, \;\; l \neq 0.  \]

Case 1. Let $l =0$.  Then ${\rm age}(l) = {\rm size}(l) =0$ and 
the untwisted 
term  equals 
\[P_{\overline{w}}^{(0)}(u,v) =     \frac{1}{uv}    \left[ \prod_{i=0}^d \frac{(uv)^{q_i}  - uv}{1 - (uv)^{q_i}}  \right]_{int}  =   \frac{1}{uv}    \left[ \prod_{i=0}^d \left( \frac{uv -1}{ (uv)^{q_i} -1}  - 1 \right) \right]_{int}  = \]
\[ =  \frac{1}{uv}  \sum_{\emptyset \subseteq J \subseteq I}   (uv -1)^{d+1-|J|}(-1)^{|J|}
 \left[  \prod_{j \not\in J} \frac{1}{(uv)^{q_j} -1}          \right]_{int}   = 
  \]
  \[ = \frac{(uv-1)^{d+1}}{uv} \left[  \prod_{i=0}^d 
  \frac{1}{(uv)^{q_i} -1}  \right]_{int}  -  \frac{(uv-1)^{d}}{uv}  
  \sum_{i=0}^d  \left[  \prod_{j\neq i}
  \frac{1}{(uv)^{q_j} -1}  \right]_{int} + \]
  \[ + \frac{1}{uv} 
  \sum_{J \subseteq I \atop |J| \geq 2}   
  (uv -1)^{d+1 - |J|}(-1)^{|J|}
 \left[  \prod_{j \not\in J} \frac{1}{(uv)^{q_j} -1}    \right]_{int}.   \]

On the other hand, 
we have 
\[  E^{(0)}_{\rm str} (X_{\overline{w}}^\vee; u,v ) = 
\sum_{J \subseteq I \atop |J| \geq 2}  \frac{(uv-1)^{|J|-1} - 
(-1)^{|J|-1}}{uv} \cdot (uv -1)^{d + 1 -|J|}  \left[
\prod_{j \not\in J} \frac{1}{(uv)^{q_i} -1} \right]_{int} =  \]
\[  \frac{(uv-1)^{d+1}}{uv} \sum_{J \subseteq I \atop |J| \geq 2} 
\left[
\prod_{j \not\in J} \frac{1}{(uv)^{q_i} -1} \right]_{int} +
 \frac{(uv-1)^{d+1 - |J|}(-1)^{|J|}}{uv} \sum_{J \subseteq I \atop |J| \geq 2} 
\left[
\prod_{j \not\in J} \frac{1}{(uv)^{q_i} -1} \right]_{int} \]

Now we use the equality $\sum_{i=0}^d q_i =1$ to obtain 
$uv= \prod_{i=0}^d (uv)^{q_i}$ and apply it in  
\[  uv  
\left[ \prod_{i =0}^d \left( \frac{1}{(uv)^{q_i} -1}  \right) 
\right]_{int}  =   \left[ uv \prod_{i =0}^d \left( \frac{1}{(uv)^{q_i} -1}  \right) 
\right]_{int}  = \left[ \prod_{i =0}^d \left( \frac{(uv)^{q_i}}{(uv)^{q_i} -1}  \right) 
\right]_{int} = \]
\[ = \left[ \prod_{i =0}^d \left( \frac{1}{(uv)^{q_i} -1} +1 \right) 
\right]_{int} =\sum_{\emptyset \subseteq J \subseteq I} 
\left[
\prod_{j \not\in J} \frac{1}{(uv)^{q_i} -1} \right]_{int}. \]
Therefore
 \[ \sum_{J \subseteq I \atop |J| \geq 2} 
\left[
\prod_{j \not\in J} \frac{1}{(uv)^{q_i} -1} \right]_{int} = (uv-1)   
\left[ \prod_{i =0}^d  \frac{1}{(uv)^{q_i} -1}   
\right]_{int}  -  \sum_{i=0}^d \left[ 
\prod_{j \neq i}  \frac{1}{(uv)^{q_j} -1}   
\right]_{int}.  \]
If we multiply the last equality by $(uv-1)^d/uv$ and add to both sides of the obtained equality the 
sum 
\[  \frac{(uv-1)^{d+1 - |J|}(-1)^{|J|}}{uv} \sum_{J \subseteq I \atop |J| \geq 2} 
\left[
\prod_{j \not\in J} \frac{1}{(uv)^{q_i} -1} \right]_{int}, \]
Then the left hand side becomes  
$E^{(0)}_{\rm str} (X_{\Delta}; u,v )$  and the right hand side becomes 
$P^{(0)}(u,v)$.
Thus, 
\[  E^{(0)}_{\rm str} (X_{\Delta}; u,v ) = P_{\overline{w}}^{(0)}(u,v). \]

\noindent
Case 2. Let $l \neq 0$.  We write $I = \{0,1,\ldots, d \}$ as disjoint
union of two subsets:
\[ I = J(l) \cup \overline{J(l)}, \; \; J(l) := \{ j \in I \; :  \; \widetilde{\theta}_j(l) \neq 0 \}, 
\;  \overline{J(l)} :=  \{ j \in I \; :  \; \widetilde{\theta}_j(l) = 0 \}. \] 
Now we have 
\[ P^{(l)}(X_w; u,v)= 
\frac{1}{uv} \left[ \prod_{j  \in \overline{J(l)}}  \frac{(uv) - (uv)^{q_j}}{(uv)^{q_j} -1}
 \right]_{int}  (-v)^{{\rm size}(l)} \left( \frac{u}{v} \right)^{{\rm age}(l)} =\]
  \[ = \frac{1}{uv} \left[ \prod_{j \in \overline{J(l)}} \left( \frac{(uv) - 1}{(uv)^{q_j} -1} -1 
 \right)
 \right]_{int} (-v)^{{\rm size}(l)} \left( \frac{u}{v} \right)^{{\rm age}(l)} = \]
 \[ \frac{1}{uv} \sum_{\overline{J} \subset  \overline{J(l)}}  (uv-1)^{|\overline{J}|}(-1)^{|\overline{J(l)}| - |\overline{J}|}  
 \left[ \prod_{j \in \overline{J}}  \frac{1}{(uv)^{q_j} -1}  
 \right]_{int} \cdot (-v)^{{\rm size}(l)} \left( \frac{u}{v} \right)^{{\rm age}(l)}  .  \]
 On the other hand,  we have 
\[  E_{\rm str}^{(l)}(X_{\overline{w}}^\vee; u,v) = \sum_{J \subseteq I \atop l \in G_J }  
\frac{ (-1)^{|J|}v^{{\rm size}(l)}}{uv} \left( \frac{u}{v} \right)^{{\rm age}(l)} \cdot (uv -1)^{d + 1 -|J|} \left[ \prod_{j \not\in J}  \frac{1}{(uv)^{q_j} -1}
 \right]_{int}.  \]
Note that $ l \in G_J \Leftrightarrow J(l) \subseteq J$ and 
$ |\overline{J(l)}| = d+1 - {\rm size}(l)$.   
 Therefore, we obtain for $ E_{\rm str}^{(l)}(X_{\overline{w}}^\vee; u,v)$ and $P^{(l)}(X_w; u,v)$ the same sums: the first sum runs over $J \supseteq  J(l)$, the second sum runs over 
the complements $\overline{J} \subseteq \overline{J(l)}$. 
Thus, we obtain 
\[ E_{\rm str}^{(l)}(X_{\overline{w}}^\vee; u,v) = 
P^{(l)}(X_w; u,v) \;\; \forall l \in \Z/w\Z. \]
If $\overline{w}$ is transverse, by \ref{mirr-trans}, this implies
\[  E_{\rm str} (X_{\overline{w}}^\vee; u, v) =  (-u)^{d-1} 
E_{\rm orb}(X_w; u^{-1},v)   \]
as it is expected in mirror symmetry. 
\hfill $\Box$

Theorem \ref{main} can be used  for constructing 
examples of Calabi-Yau varieties $X$  which do not have 
mirrors:  

\begin{exam} Take the weight vector  with $IP$-property
$\overline{w} = (1,1,2, 4,5)$. Then $w = 13$. We note 
that this weight vector is not transverse. By Theorem \ref{main}, 
we can compute the stringy $E$-function of a Calabi-Yau 
compactification $X_{\overline{w}}^*$ of the affine hypersurface 
$Z_{\overline{w}} \subset (\C^*)^5$ defined by the Laurent polynomial 
\[ f_{\overline{w}}^0({\bf t}) = \frac{1}{t_1 t_2^2 t_3^{4}t_4^5} +
t_1 + t_2 + t_3 + t_4.   \] 
All nonzero elements $l \in \Z/13\Z$ have size $5$ and 
determine twelve terms that give rise to polynomial $ -v^3 - 5uv^2  - 5u^2v - u^3$. 
The untwisted term ($l=0$) 
\[  \left[  \left( \frac{1-(uv)^{\frac{12}{13}}}{1- (uv)^{\frac{1}{13}}} \right)^2 \cdot   \left( \frac{1-(uv)^{\frac{11}{13}}}{1- (uv)^{\frac{2}{13}}} \right) \cdot  \left( \frac{1-(uv)^{\frac{9}{13}}}{1- (uv)^{\frac{4}{13}}} \right) \cdot  \left( \frac{1-(uv)^{\frac{8}{13}}}{1- (uv)^{\frac{5}{13}}} \right)
\right]_{int}
\]   
is not a polynomial, because its value for $u=v=1$ is 
$\frac{1092}{5} \not\in \Z$ (cf. \cite[Example 1.13]{BS20}). This 
is exlained by existence of terminal singularities on minimal models of $Z_{\overline{w}}$ \cite{DR01}.  
Therefore, $E_{\rm str}(X_{\overline{w}}^*; u,v)$ is not 
a polynomial, and $X_{\overline{w}}^*$ has no mirror. 
\end{exam}

\end{document}